# The Machine as Data:
## A Computational View of Emergence and Definability

S. Barry Cooper[1]

School of Mathematics, University of Leeds,
Leeds LS2 9JT, UK

*In 1936 Turing developed the definitive theory of universal classical computers. His motivation was not to build such a computer, but only to use the theory abstractly to study the nature of mathematical proof. And when the first universal computers were built, a few years later, it was, again, not out of any special intention to implement universality. They were built in Britain and the United States during the Second World War for specific wartime applications. The British computers, named Colossus (in which Turing was involved), were used for code-breaking; the American one, ENIAC, was designed to solve the equations needed for aiming large guns. The technology used in both was electronic vacuum tubes, which acted like relays but about a hundred times as fast. At the same time, in Germany, the engineer Conrad Zuse was building a programmable calculator out of relays — just as Babbage should have done. All three of these devices had the technological features necessary to be a universal computer, but none of them was quite configured for this. In the event, the Colossus machines never did anything but code breaking, and most were dismantled after the war. Zuse's machine was destroyed by Allied bombing. But ENIAC *was* allowed to jump to universality: after the war it was put to diverse uses for which it had never been designed, such as weather forecasting and the hydrogen-bomb project.*

      — David Deutsch: *The Beginning of Infinity*, Allen Lane/Penguin, London, 2011, p.139.

## 1. Language, Representability and Universality

Turing's (1936) paper on computable numbers has played its role in underpinning different perspectives on the world of information. On the one hand, it encourages a digital ontology, with a perceived flatness of computational structure comprehensively hosting causality at the physical level and beyond. On the other (the main point of Turing's paper), it can give an insight into the way in which higher order information arises and leads to loss of computational control — while demonstrating how the control can be re-established, in special circumstances, via suitable *type* reductions (of which more later). Our preliminary aim is to examine the classical computational framework more closely than is usual, drawing out lessons for the wider application of information-theoretical approaches to characterising the real world. The problem which arises across a range of contexts is the characterising of the balance of power between the complexity of informational structure — with emergence, chaos, randomness and 'big data' prominently on the scene — and the means available (simulation, codes, statistical sampling, human intuition, semantic constructs) to bring this information back into the computational fold.

---

[1] This article is based on an invited talk by the author at the *First International Conference on Logic and Relativity: Honoring István Németi's 70th Birthday,* held at the Alfréd Rényi Institute of Mathematics, Hungarian Academy of Sciences, September 8—12, 2012. Research and preparation of this article was supported by a John Templeton Foundation research grant: *Mind, Mechanism and Mathematics*, July 2012 — August 2015. The author is grateful to the two anonymous referees for their detailed questions and suggestions which led to a number of improvements in the final version of the article.

First impression is that Turing's universal machine is not a machine at all, rather a schematic computer program. The history, decisively moulded by the intervention of John Von Neumann with his so-called 'Von Neumann architecture', points to the subtlety of the interaction of concepts. In reality, machine and program are hard to disentangle. We observe a confluence of logic and matter in which all computational elements (including program) are embodied, but in which the structure of the embodiment, guided by the logic, is hugely important. The importance is a practical one of feasibility, where the logical structure of the universal machine can be arrived at in different ways. The routes of Von Neumann, or less influentially that of Turing's ACE report, were the ones that held the key to the future. Others, involving physical fragmentation via punched cards or hand-sorted plug boards were doomed, from the start, to obsolescence.

The David Deutsch quotation above is chosen to illustrate how quite brilliant scientists with unimpeachable academic credentials can underestimate the subtlety of what is going on, and the role of the logic. The universal machine was not just proficient with language enough to be programmable. It was not just equipped with a large enough range of actions, and accepting enough of instructions to perform those actions, to be able to carry out any properly composed program devised. Turing's universal machine contained a *dictionary* of machines! One could carry on a conversation with it about programs and what they deliver. Anyone who has been sent a document unreadable on her computer due to a missing piece of software or absent macro will understand the value of the stored program computer. And go back to the days before the internet, where a physical journey to a store, or reliance on the post office, was needed to fill the deficiency, and one has a strong sense of the value of embodying the software in-house, rather than via disconnected embodiment, effectively speaking a different language.

Hence the 'machine as data' paradigm which the universal machine delivered, but which is both informative and misleading. The Turing trick was to devise a structure which one understood to handle instructions, and what these acted on, in just the same way. The program for the universal machine was a sophisticated one, and is embodied in a modern-day computer in a sophisticated way, embodied 'in-house', certainly, but using a complex balance between architecture, system and programming. All is embodied. There is no such thing as separate data, anymore than one can separate mind and body. The key problem for philosophy of mind is

*understanding* the connection. That for computer science and artificial intelligence is *managing* it.

All is dependent on both logical structure and creative engineering. A 2014 article by Thomas Haigh in the *Communications of the ACM* is entitled *Actually, Turing Did Not Invent the Computer*. Behind this truism one can discover another aspect of the same misunderstanding apparent in David Deutsch's impression of the versatile ENIAC 'jumping' to universality via its programmability. The mathematical model is far from comprehensive and unique in its prescriptive detail. And the ad hoc ingenuity of the engineering is far from independent of the logic. Punched cards and plugboards do not disqualify one from programming anything a universal machine can do. According to the terminology, machines such as Colossus and others mentioned by David Deutsch may well be *Turing complete*. One may even, modulo some suitable interpretation of the physical structure of the machine, inductively *prove* that. But universality needs to be built in. It does not happen by accident. And even if it did, it would be hard to recognise without knowing the underlying design considerations.

So far, we have examined what is essentially a reduction of machine to data, observed the convergence of the two, and seen how powerful the outcome might be. We have seen the importance of the language in the programming. The significance of representability of the machine via its reduction to simply definable logical structure. And the universality that emerges from this so powerfully. And for some, this might seem the end of the story.

One of Turing's special interests was type theory. This concerns the same typing of information that Bertrand Russell resorted to, in order to bring to the abstractions of the higher mathematics the same sort of consistency that we observe in the real world. As Gödel summarises in his *Russell's mathematical logic* in 1944:

> *By the theory of simple types I mean the doctrine which says that the objects of thought … are divided into types, namely: individuals, properties of individuals, relations between individuals, properties of such relations, etc. … , and that sentences of the form: " a has the property φ ", " b bears the relation R to c ", etc. are meaningless, if a, b, c, R, φ are not of types fitting together. Mixed types (such as classes containing individuals and classes as elements) and therefore also transfinite types (such as the class of all classes of finite types) are excluded. That the theory of simple types suffices for avoiding also the epistemological paradoxes is shown by a closer analysis of these.*

One can see the typing playing out in a very specific context in the 1936 paper. As well as the reduction of machines to data via coding of logical structure, one sees correlation of computations, successful and otherwise, into an incomputable type-one object — what we now read off as an incomputable real number derived from the so-called (unsolvable) *halting problem* for the now standard version of a universal Turing machine. By *incomputable* (an old word with a 20th century interpretation) we mean that the halting real, viewed in a standard way as a type-1 function, cannot be computed by any Turing machine. The discovery of such an object might have had more impact in the years after 1936 if it had been something more clearly observable in our everyday lives. It really needed a natural embodiment in nature. Unfortunately, it was the type-ascending lever to incomputability that seemed to belong essentially to a niche mathematical world of *recursive function theory*, one not even of huge interest to logicians themselves. In the new world so successfully ruled by the Turing computational model, there was not even the motivation to look around for the embodied halting real. There was plenty could not be computed, but one could ascribe that to common or garden complexity.

The mathematics giving the incomputability of the halting real is impeccable. It takes us surely and undeniably from solidly embodied machine, familiar in theory and in everyday experience, to an object which cannot be computed in detail by its physical host. It is an abstract object, We cannot see it. It has no clear relevance to anything in the real world. Import the halting problem into the real world and the unbounded time which is the root of the incomputability evaporates in the sunshine of a familiarly local context. Even infinity — despite the title of David Deutsch's book — is of doubtful relevance to the world we experience. But the sense of mismatch between the digital and the observable wider reality persists.

## 2. Relevance, Relativity and Digital Ontology

So we have a paradoxical situation. On the one hand, the mathematics of Turing's 1936 paper dominates the modern world via the stored program computer and its accompanying culture. Algorithmic approaches to everything, appropriate and inappropriate, are in vogue. What is different to previous cultures, ruled by religion or ideology, is that while the content of rules may be less micro-managed in more liberal western societies, there is an increasingly ingrained faith in the efficacy of algorithmic controls and procedures per se. Whole layers of administration have grown up around this outlook. At the same time the mathematics underlying the key conclusion of the 1936 paper — that not everything is computable — is pruned off, along with much else, above the type-0 level.

On the other hand, there is a growing groundswell of dissatisfaction with the pervasive flatness of this implicit digital ontology. The unsolvability of the Halting Problem may be relegatable to the playground of recursion theorists; but the cultural and scientific weight of evidence impelling recognition of an associated higher type ontology — driving much of Turing's later work — is becoming hard to ignore. Our growing understanding of the profound consequences of the young Turing's simple diagonalisation of the partial computable functions is almost complete.

This dissatisfaction with digital ontology comes from those concerned with creativity and human individuality of thought; from academic disciplines, such as economics, in which the mathematics has often failed and is suspect; from social scientists facing the reality of a hybrid discipline; and from those, like Turing in the 1940s, who have accepted the need for computers to be informed by interaction with human common sense. It also arises more formally via philosophical inspection of the deficiencies of digital ontology. In Luciano Floridi's (2011) *The Philosophy of Information* he summarises (chapter 15):

> *… a digital (or analogue) ontology is not a satisfactory approach to the description of the environment in which informational organisms like us are embedded.*

The paradox we have in mind is that on the one hand mathematicians and computer scientists, and scientists more generally, tend to peripheralise the mathematics of incomputability in favour of a restrictive digital ontology. While those who might be expected to find in the mathematics beyond that of the computable some support for their broader view of the world, bypass the full richness of structure available, unfamiliar with the technical detail, and repelled by what they see as the dominance of the classical paradigm. It is somewhat sad to see the world beyond the computable largely relegated to the status of our fairies at the bottom of the garden. It may be that mathematical modelling is irrelevant to the specifics of understanding or predicting the behaviour of, say, people. But basic fundamentals can determine important global characteristics, such as the nature of the limitations on our knowledge. These comments set out our approach. A number of the issues we touch on here are discussed in depth in rather different terms in the extensive literature. We try to avoid the technical terms from philosophy that might divert attention from the aim, which is to bring some less than familiar mathematics to play on fundamental questions of structure and causality in the real world.

It is true that we have no firm grasp on 'reality' (whatever that is). For some, reality is so elusive, that for them all that exist are are narratives whereby we interpret and survive in the world. We take a broadly scientific perspective, informed by what science does and how it does it. We recognise the bounds of the representational resources available to a science that deals in clear descriptions of reality that can be tested and shared. At the most basic level, these are essentially computational and observational bounds. Here, with the above description of types from Kurt Gödel still in mind, we recognise and name individuals, mainly at the type-0 or type-1 level, only occasionally straying beyond type-2 — covering manifolds, vector fields, etc. — in our discussions.

In practice, much is approximated and dealt with computationally at the familiar type-0 level. And in the field, this is often achieved through computationally determined sampling of phenomena, ones which potentially might require a higher type description. Richard Feynman memorably pronounced in his 'keynote' talk at the 1981 conference at MIT on the `Physics of Computation':

> *It is really true, somehow, that the physical world is representable in a discretized way, and … we are going to have to change the laws of physics.*

Although this might reassure us about the usefulness of the classical computational model (and the related quantum computing version), experience and — as we will see later — such type reduction, can result in a faulty guide to the reality. Also, representation is inevitably *relative*, to axes or observational base. We look at the computational consequences of this later, when we have established an extended computability theoretic framework more suited to the informational complexity of our universe.

It is in the social dimension that extending the Church-Turing thesis meets specially informative challenges. Is the social context, which is so basic to our everyday concerns, to be relegated to a non-scientific play area beyond serious scientific attention? What can we extract from the following example of John Searle (in his 2010 *Making the Social World: The Structure of Human Civilization*, pp.88-89)?

> *Let us suppose in a pub that I get up from the table and go to the bar and order three beers. I then carry the beers back to the table and set them down. I say, "This one is Sally's; this one is Marianne's; and this one is mine." Now this would not appear to be a very remarkable metaphysical effort, but in fact it has remarkable properties. By making these utterances, I have in fact created new rights. Indeed, they were more than statements; they were deontic commitments of a very special kind. I created a reality according to which, for example, Sally has certain rights that Marianne does not have and Marianne has rights that Sally does not have. This would come out in the fact that if Marianne tried to drink Sally's beer, Sally would have a legitimate complaint.*

What is clear from this apparently simple example is that something beyond the simplicities we are familiar with is going on. The imposition of an agreed outcome is dependent on a non-local appreciation of the context within which this simple ritual is taking place. Is the phenomenon beyond the reach of analysis, relegated to the company of our fairies at the bottom of the garden? Is this the only alternative to the unfeasibility of digital computation? The sense is that there is a familiar semantical content to the situation that is real in its operative relevance, but demands higher level analysis for a better understanding. And we do not expect the Turing machine, unaided by human intelligence, to be much help with this. What we are seeing here is what Searle tells us is a *Status Function Declaration*. For us, it is the calling up of what is clearly a shared 'computation' of an outcome: with the task in hand that of a convincing removing of the quotation marks. The relatively simple outcome described by Searle can be viewed as the application of an 'analogy' in the sense of Hofstadter and Sander in their *Surfaces and Essences — Analogy as the Fuel and Fire of Thinking*, from 2013. It is the identification and packaging of a body of relatively complex data in a form suitable for the application of the words committing us to the outcome.

Despite finding Searle's example of the beers interesting and somewhat suggestive, the reader might not find it entirely convincing as a case for taking seriously the existence of higher levels of information that are not readily reducible to a digital level. The reader used to *computationally* scheduling complex daily activities might identify quite similar things in the world of computer networks — e.g. a server might grant to clients the rights to access certain resources by means of individualised 'tickets' of some kind — and no one thinks there is anything very mysterious to be explained about this. The speculation — essentially a 'proof by example' — might be that open systems can be reduced very generally to data algorithmically organised. But as Floridi remarks:

> *. . . a digital (or analogue) ontology is not a satisfactory approach to the description of the environment in which informational organisms like us are embedded.*

As Floridi describes, it is the interface between 'inforgs like us' and the structures we seek to impose that make necessary analysis in terms of the *Levels of Abstraction* (LoA's), and which fragments reductive thinking. Here he is describing in *The Philosophy of Information* (pp.347–348) the second of "four advantages" of the method of levels of abstraction:

> *... specifying the LoA means clarifying from the outset, the range of questions that (a) can be meaningfully asked and (b) are answerable in principle. We saw that one might think of the input of an LoA as consisting of the system under analysis, comprising a set of data; its output is a model of the system ...*

While (p.349):

> *... by accepting an LoA a theory commits itself to the existence of certain types of objects, the types constituting the LoA ...*

The informational infrastructure we bring to the computational context changes both incident and rules for dealing with such particularities. Searle's *Status Function Declaration* is not presented as a barrier so much as an invitation to creative thinking concerning the moving social landscape. It invites us to open out the situation into one where in his example of beers there is confusion over description and ordering of items, changing of states of mind and incidents involving spillage, or remembrance of an allergy, with consequent negotiation or dispute over what was ordered. What takes the situation beyond the simplicity of the algorithm is not the lack of definition of the algorithm, nor a carelessness in the summoning up of universality, but the invocation of the higher order analysis arising from the eroding of the closed system so basic to the algorithmic approach, via the relational potential of the human inforgs. The attempt to incrementally eliminate the potency of Searle's example via a transfer of say the ticket algorithm has the flavour of Turing's 'failed' 1939 attempt to replace 'intuition' by 'ingenuity' in his iteration of Gödel's incompleteness into the transfinite. One finds a more literary take on the relationship between different 'Levels of Abstraction' in D. H. Lawrence's 1926 novel *The Plumed Serpent:*

> *Let us seek life where it is to be found. And, having found it, life will solve the problems. But every time we deny the living life, in order to solve a problem, we cause ten problems to spring up where was one before. Solving the problems of the people, we lose the people in a poisonous forest of problems.*

Also significant: If one tried to digitise the higher order social computation implicit in Searle's example, one would lose the flexibility and more general applicability, while introducing an unrealistic expectation of infallibility. As Turing said in a talk to the London Mathematical Society in 1947 (see Hodges *Alan Turing: The Enigma*, p.361):

> *... if a machine is expected to be infallible, it cannot also be intelligent. There are several theorems which say almost exactly that.*

The role of the data is not simple, and the shared nature of the identification of its role is fragile, subject instabilities flowing from possible observational disparities between those acting out Searle's scenario. However, there are quite definite rules at work. These correspond to the sort of rules underlying our survival in the physical world, applying in the more complex social world that Searle is discussing. Can we find a mathematics which mirrors the computational aspect along with the less secure handling of the 'big data'?

Interestingly the English legal system, based on a complex combination of legislation and carefully formulated judicial precedent, contains within it an implicit recognition of the limitations of a purely algorithmic approach. English common law is essentially created by judges sitting in court and applying (with ad hoc refinement and augmentation) legal precedent (*stare decisis*, or 'to stand by things decided') to the facts before them. The key to its effectiveness is the role the process gives to decision based on human judgement rather than purely written legislation. We observe all the hallmarks of the institutionalisation of an appeal to higher order computation, encapsulated in Hofstadter-Sander use of analogy. The process has been developed so as to avoid the ad hoc on the one hand, and the arbitrary on the other, via a computational process in which it is the identification of the appropriate data which legitimises the computation. The connection with the Searle example with its engagement with higher order information should be clear by now.

A persuasive view of the limitations governing the representational framework is that the resulting descriptions do *define* constraints on reality, but in a potentially incomplete way. This is seen in its most convincing form in physics. Without some very basic mathematics, we are open to the argument that there may be many different narratives, which will fit observation from different points in history. This could lead to a view of Kuhnian paradigms not just *changing* under the influence of conflicting social pressures, but entirely *based on* such social dynamics. This would be a view recognizable from Richard Rorty's 1988 essay on *Science as Solidarity* (revised and reprinted in *Objectivity, Relativism, and Truth*, Cambridge University Press, 1991) where he says (p.39):

> *We cannot, I think, imagine a moment at which the human race could settle back and say, "Well, now that we've finally arrived at the Truth we can relax." …*
>
> *Pragmatists would like to replace the desire for objectivity — the desire to be in touch with a reality which is more than some community with which we identify ourselves — with the desire for solidarity with that community. They think that the habits of relying on persuasion rather than force, of respect for the opinions of colleagues, of curiosity and eagerness for new data and ideas, are the only virtues which scientists have. They do not think that there is an intellectual virtue called "rationality" over and above these moral virtues. …*
>
> *My rejection of traditional notions of rationality can be summed up by saying that the only sense in which science is exemplary is that it is a model of human solidarity.*

There is a route out of this (what may seem to many) counter-intuitive view of science. Maybe we can identify *something* that we can accept as basic to the representational framework.

# 3. The Computational Structure and Embodiment of Higher Order Information

Our intuitive finitism is countered by the algorithmic underpinning of so much causal structure in our world, the computability of which helps the world make sense and provide the basic ingredients for our survival in it. What we then become subject to is the infinitary mathematics that the algorithms drag with them. When the interactions of (sub)atomic particles are governed by forces ranging from the weak and strong nuclear to the gravitational, one is no longer dealing with the mathematical simplicity of a finite or even discrete relational structure. The mathematical costume within which the universe dances before us has a logical and mathematical complexity supporting layer upon semantical layer. Typically, fluid dynamics is equipped with computable rules governing form via appropriate differential equations. The reality we observe from the sea shore — or looking out into space — fills us with respect for the mathematical descriptions underlying the incipient unpredictability.

We live in the embrace of mathematical structures called up by the algorithms, such as are so effectively described by programs like those Turing packed into his 1936 article, and which start by delivering the richness of the computable reals underlying science. Open the Pandora's box of infinity and its mathematics, and the result is that the typing of information, so basic to Bertrand Russell's restoration of the integrity of the way we deal with information and the formality of rationality, becomes an unavoidable reality. This is all very well for the more abstract-minded amongst us. For those of us who look for real world significance, we ask where can we actually see this hypothesised type structure embodied for us all to see? What is there to convince us that there is anything here of relevance for the realist?

Entering the everyday world characteristically free of mathematics, but not of rigorous thought, we could do worse than go to Luciano Floridi's excellent book on *The Philosophy of Information*. Here we have a refreshing dispensing with the digital ontology imposed on us by mathematicians, computer scientists, and physicists of a similar frame of mind to David Deutsch. A source of misgivings surrounding the digital ontology, from the logician's point of view, is the familiar experience that if something occurs mathematically, then it is just a matter of time before it crops up in the material world in some guise. From the point of view of the philosopher with a strong sense of what is real, the flatness of the digital ontology is just too out of kilter with what we so powerfully perceive, in the words of Wittgenstein, "to be the case". This intuition is tied to an ingrained faith in the semantical power of language, with all its nuances and resistance, both heuristically and theoretically, to being handled by computers.

Though the reach of mathematical language in the world we encounter in everyday life may be limited by what we can computationally verify, such considerations are not so relevant to the use of philosophical argument. The validity of a statement by Wittgenstein or Heidegger or Derrida is found to have validity in a quite different way to that of say Newton or Maxwell or Darwin. But — there there is a well-worn mathematical route between algorithmic content and natural language, and back. This is encountered schematically in the various logical characterisations of hierarchies of statements in terms of their computability-theoretic underpinnings. At a basic level the best-known example is Post's Theorem. On the other hand, the fragmentation of the world of the working scientist may break the linguistic grasp of transitions between different levels of scientific knowledge, but there remains a strong sense of causal connection between (say) the quantum and the biological levels of our world. Despite the impracticality of algorithmically capturing the translation between levels in particular cases, there is little reason to doubt that natural language is as rooted in the basic structure of the universe as is mentality supervenient on the physical characteristics of the brain.

The aim above is to restore some consistency between the world of language and its semantical accretions; and between the equally sophisticated but (in terms of its carelessly abstract form) brutal mathematics based on, but not contained by, the science. An obstacle to the philosophical acceptance of the mathematics — besides an instinctive distrust of the mathematician's tendency to over-simplification and invalid reductions — is an apparent informational similarity of type at different levels of abstraction; and the sheer complication of structure, apparently far more complex than the 'levelism' (such as disowned by Floridi) implicit in mathematical hierarchies. The logical structure maybe of limited practical value in transporting meaning between levels in the sense of Floridi, but is invaluable as a means of introducing inductive structure to a computationally complex environment.

To summarise — a mathematical response to such reservations would point to the value of an analysis of hierarchy as a basis for the deconstruction of more complex non-linear relational frameworks. And to specially prioritise the task of understanding the character of the relationships between LoAs and the frequently encountered qualitative transitions demarcating the informational terrain. We need to raise our awareness of the mechanisms for type reduction, and for informational 'phase transitions' which certainly cross semantic/computational barriers in ways which do not depart our observational domain. Out of the looked for reconciliation of outlook and objective can develop a coherent approach to *dealing with* levels of abstraction — particularly those whose relationships to the wider context are imperfectly characterised, such as the quantum level, mentality, natural laws and higher-order structure in the universe, and emergent phenomena more generally. The understanding one is looking for at the informational and epistemological level, as always, is informed by embodied analogues. David Ruelle, an originator of the term *strange attractor*, comments in his 2007 book on *The Mathematician's Brain* (pp. 120-121):

> *Considerable efforts are currently under way to obtain a "theory of everything," allowing, in principle, to understand all observed physical phenomena. When such a theory has been obtained, it will be possible to compute every physical quantity … It would then seem that the most interesting part of physics will be over, the rest being "just calculations." But this is not the case, because there are important conceptual problems in physics that go far beyond the discovery of fundamental laws …*
> 
> *Think now of understanding the properties of water, given that you know the fundamental laws of the mechanics of water molecules. You would, for instance, like to understand phase transitions: why, when you change its temperature, does water suddenly freeze to ice or boil to vapor? You would like to compute the viscosity of water (its resistance to deformation), and you would like to understand turbulence. … The properties just mentioned are emergent properties. They are not properties of one water molecule or ten water molecules — they appear in the limit of infinitely many molecules. It is true that in the lab you always work with a finite amount of water, but the number of molecules in a liter of liquid is huge, and the properties of interest are (in first approximation) those of an infinite system.*

In 1939, Turing (for whom Ruelle expresses a high regard in his book, see pp. 83-84) set out to scale the epistemological barrier presented by Kurt Gödel's incompleteness theorem. In doing so, he traced a route through a (mathematical) phase transition which, in essentials, contained the ingredients of other such *observer-related* traverses in superficially different contexts. The machinery was computational in detail, pinned to computable (or what Kleene called *recursive*) ordinals, but put together in an extended relationship that would vie with the familiar halting real for a similar level of incomputability.

A key element in subsequent deconstructions of natural 'phase transitions' at higher type was the one-page of this long 1939 article devoted to the introduction of the oracle Turing machine. In terms of the science, generally framed in relation to the real numbers, this was exactly what was needed to model the computable causality basic to most of what we know of the natural world. And the rich theory that it gave rise to would provide basic levers to the modelling of phenomena which do not fit in this classical setting. What distinguishes the oracle Turing machine from the classic 1936 incarnation is the recognition that we frequently compute *relative* to information we may derive from more-or-less trusted sources, without detailed knowledge, computational or otherwise, of its origins. And what makes it unique, giving rise to a unmistakable modelling of scientific relationships, is its capturing of the intrinsic character of scientific computation of causal relationships, complete with its commonly observed basic continuity and finite approximability. There are other models mathematically based on careful observation of causal relationships specific to particular contexts — such as Rafael Sorkin's causal sets based on the physics — but the very specificity of the scope of success of such models betrays their lack of fundamentality.

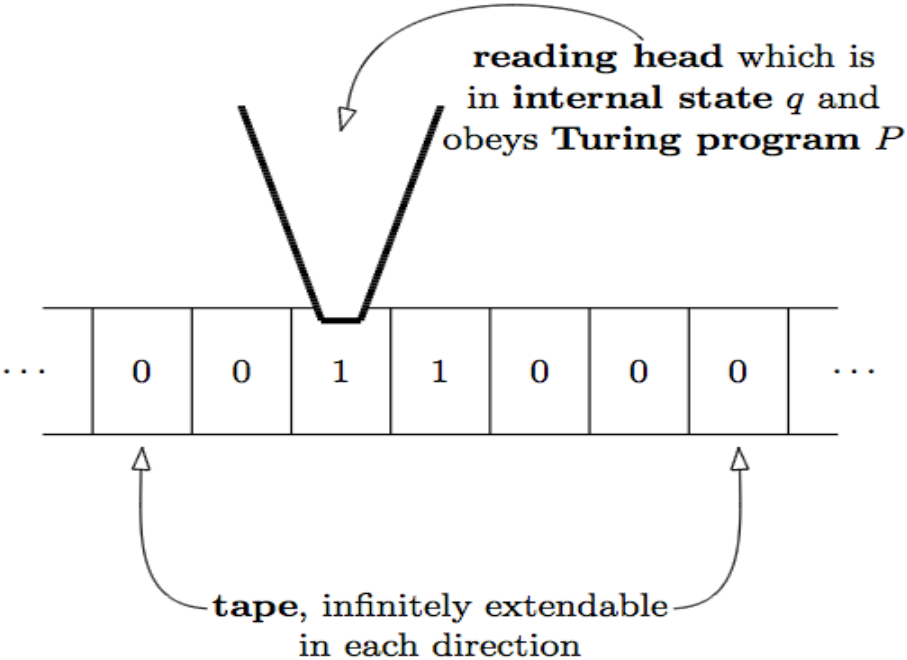

Intuitively, the computations performable by an oracle machine are no different to those of a standard Turing machine — except that auxiliary information can be asked for, via questions about increasingly long finite parts of a given real number (termed the *oracle*). A typical machine using a binary language, with a subdivided extendable tape and a reading-writing head, is shown above. The program for the machine consists of atomic instructions asking the reading head to identify what it is reading, and its current internal state — *of mind*, to follow Turing's mental model — and move left or right one cell of the tape, or print or erase a tape symbol. To extend such a finite program for a regular Turing machine to one with the oracle consulting facility, one merely allows a new type of instruction which says "count up the number *n* of 1's currently printed on the tape, and ask if entry number *n* of the oracle real is a 1 or not — and go into the state named in the instruction as corresponding to the binary possibility pertaining".

This model is seen most clearly in today's digital data gathering, whereby one is limited to receiving data which can be expressed, and transmitted to others, as information essentially finite in form. But with the uniformly computational character of the model comes the capacity to collate data in such a way as enable us to deal with arbitrarily close approximations to infinitary inputs and hence outputs. This can give us a computational counterpart to the computing scientist working with real-world data. If the different number inputs to the oracle machine result in 0-1 outputs from the corresponding Turing computations, one can collate the outputs to get a binary real computed from the oracle real, the latter now viewed as an input. This gives a partial computable (p.c.) functional $\Phi$, say, from reals to reals. One can obtain a standard list of all such functionals, where as usual, there is no way of knowing in general whether a given computation ever halts, or, for that matter, provides the sort of information you wanted. We say the output real is *Turing computable* from the real input to the functional $\Phi$.

The set of all such p.c. functionals $\Phi$ acting on the set $\mathbb{R}$ of real numbers gives the *Turing Universe* as a schematic model of the basic causal structure of the real world. The American mathematician Emil Post (1948) gathered together into equivalence classes those sets of binary reals which are mutually Turing computable from each other. There was an ordering on the equivalence classes induced by the original relation of "Turing computable from". Mathematically, this delivered a much tidier structure, the upper semi-lattice of the *degrees of unsolvability*, which subsequently became familiar amongst logicians, and just a little beyond, as the *Turing degrees*.

There are obvious parallels between the Turing universe and the scientific view of the material world. Most basic, scientific theories describe the world in terms of real numbers. This is not always immediately apparent. Nevertheless, scientific theories consist, in their essentials, of postulated relations upon reals. These reals are abstractions, and do not come necessarily with any recognisable metric. They are used because they are the most advanced presentational device we can practically work with, although there is no faith that reality itself consists of information presented in terms of reals. Also, the relativity of the informational image of the world has to be factored in. There is in general no absolute picture of a large informational structure. Take an observer situated at the zero of a representation of euclidean 1-space, and the computability or otherwise of other locations on the line may radically differ in specificity from those points an observer at another point reaches out for. Of course, the plethora of automorphisms of the real line means that the actual structure of the informational world of one observer will be isomorphic to that of another. But with more complexity of structure this may not be the case anymore. This is what one would expect from a familiarity with physics of relativity.

In general, physicists would not allow the computational content of their activities a role in dictating a physical embodiment the real numbers, since any physical quantity can only be measured to a finite precision, and can only be measured to a *low* finite precision (a handful of significant figures). What one is referring to here is the mathematical *model* that *hosts* the scientific activity. Theoretical models pre-bounded in terms of resources described (as opposed to specifically calculated) would present unacceptable difficulties. It is not that '√2' or 'e' systematically approximated would inflate errors — it is that the destruction of mathematical structure via the replacement of logically explicit descriptions of objects by numerical approximations would carry through to a corresponding destruction of logical content of outcomes, leading to a potential catastrophic degeneration of the whole scientific enterprise. This takes on significance when one has descriptions of scientific parameters which do not immediately lead to computable approximations. But here we enter territory where there is no consensus around the role of the computability theory in the scientific context. Easier just to watch the mathematical physicist at work, and observe the extent to which the computer is used to do more than correspond with colleagues and format the paper written at the end of the day. Relative Turing computation reflects the real context in basic detail, and is sensitive in its formulation to the role of approximation in the mathematics called upon.

Globally, there are still many questions concerning the extent to which one can extend the scientific perspective to a comprehensive presentation of the universe in terms of reals—the latter being just what we need to do in order to model the immanent emergence of constants and natural laws from an entire universe. Of course, there are many examples of presentations entailed by scientific models of particular aspects of the real world. But given the fragmentation of science, it is fairly clear that less natural presentations may well have an explanatory role, despite their lack of a role in practical computation.

The natural laws we observe are largely based on algorithmic relations between reals. Newtonian laws of motion will computably predict, under reasonable assumptions, the state of two point masses moving under gravity over different moments in time. And, as previously noted, the character of the computation involved can be represented as a Turing functional over the reals representing different time-related two-particle states. There is no claim that the model is more than successful at a certain level, any more than the specificities of the Newtonian laws themselves are. One has no idea what reality actually is, one merely has a level at which one occupies it, and an epistemological toolbox of some limited practicality and reliability. We stress again (following Floridi and other philosophers), that there can be no assumption in this discussion that reality *is* real numbers.

Further, one can point to physical transitions which are not obviously algorithmic, but these will usually be composite processes, in which the underlying physical principles are understood, but the mathematics of their workings outstrip available analytical techniques. Such mathematical phase transitions may need very little in the way of composition of ingredients. For instance, back in 1970 Georg Kreisel (in a footnote to his article on *Church's Thesis: a kind of reducibility axiom for constructive mathematics*) suggested a collision problem related to the 3-body problem, which might result in "an analog computation of a non-recursive function". For some, the dependence here on the embodiment of an analog device of unbounded precision may weaken the challenge to a computable universe. The computability theorist will stress the basis of the mathematics of the incomputable in computable approximations, which impinge far more surely on our physical and mental world. Kreisel knows this, and engaged viscerally over a long period with the physicality of abstract computation, and his thinking here is worth unpacking from the conciseness of his footnote. His take on the physical summoning up of infinitary mathematics, and its incipient incomputability, has been echoed in more recent times by Donald Saari and Jeff Xia in an interesting take on the Painlevé conjecture in their May 1995 article on *Off to Infinity in Finite Time*.

What is important about the Turing universe is that it has a rich structure of an apparent complexity to parallel that of the real universe. At one time it was conjectured (the so-called *bi-interpretability conjecture*) that the definable relations on it had a mathematically simple characterisation closely relating them to those of second-order arithmetic. Nowadays, with the passage of years, it appears more likely that the Turing universe supports an interesting automorphism group, echoing the fundamental symmetries emerging from theoretical descriptions of the physical universe. On the other hand, there are rigid substructures of the Turing universe (see Slaman's *Degree Structures* paper from 1990-91) reminiscent of the classical reality characteristic of our everyday world. And the definable relations hosted by substructures provide the basis for a hierarchical development. As we noticed in the simple example of the information content of the real line — a view consistent with Floridi's ISR "*informational structural realism*, a version of OSR supporting the ontological commitment to a view of the world as the totality of informational objects dynamically interacting with each other" — we usually discuss information defined or (more particularly) computed from information of a more basic character. By the 'hierarchical development' above we mean via the application of Turing definability via a hierarchical iteration of representations limited to those acceptable in a scientific context. This is not to deny the 'reality' of higher order information beyond what we can even sample, or to deny the 'downward causality' embodied by the language and its computability theoretic underpinnings.

Our overview, which one might characterise as science with consciousness, is not intended to distract credit from the hugely important enterprise of drawing out the specific descriptions relating to both large and small LoAs hosted by our universe. An important part of this focus on key aspects of the world we live in is bound to be some real validation of our take on the computational content of the real universe. One of the more exciting and challenging aspects of this, bringing us face-to-face with computationally novel physical dynamics, is the work on the consequences and nature of relativity pursued by István Németi and the other members of the Budapest group led by him and Hajnal Andréka. Amongst many important articles, we would single out the 2006 paper of Németi and Dávid on *Relativistic computers and the Turing barrier*, and, from Hajnal Andréka, Judit X. Madarász, István Németi and

Gergely Székely, *A logic road from special relativity to general relativity* (*Synthese*, 2012).

It is worth noting that although we identify a key role for the typing of information in accordance with the observational distinctions between local and global dynamics within our universe, the epistemological complexity and fragmentation we encounter is not characterised purely in terms of type. Observation itself entails a physical relationship whose embodiment will take place at a particular level of *computational* relatedness. Our observational activity does itself embody considerable, and very effective, extraction of emergent information and the emergent relations on it in qualitatively comparable form to that which it is emergent from. In our everyday world, we may retain distinctions between physical observables and their semantical content. More often, such distinctions disappear as we move continually between different LoAs.

The termite cathedral may exhibit form which is clear to us as observers and not to the termite. But a purely digital representation of the termite cathedral — a photograph say — becomes an object no different in character to that of any other digitally represented entity. We have a semantical appreciation of the overall form of the cathedral via the large representational resources of the brain; but this appreciation takes place at a higher level, at which the brain may mimic the type reductional facility of the camera via its evolved cognitive embodiment.

By extension, biology may be emergent from quantum mechanics via some route structure, but qualitatively it is a science permeated by reasonably effective digital informational structures, like any other. There are even interactions discernible between the quantum and biological worlds (see the 2009 Arndt, Juffmann, Vedral article on *Quantum physics meets biology,* or Rieper, Anders, and Vedral in 2010 on *Entanglement at the quantum phase transition in a harmonic lattice*). But the barrier one encounters to translation of one world into another is computability theoretic, arising from a type transition or more, but ones not always explicitly recognisable. This is a feature of examples of phase transitions cited by David Ruelle. The world does deceive us with its apparent informational flatness. A quantum complexity can give rise to a Linnean orderliness with a semantics not immediately indicative of the underlying informational relationships. Normality may be built on something far more mysterious and inaccessible — and more dangerous — than we wish to contemplate. The work of Vedral and his co-workers is just one small pointer to the incidence of causal bleeding between relatively well understood scientific levels. It may not be "turtles all the way down", but there is an undeniable possibility of embodied interactive informational levels beyond the relative instantiation of the human organism.

# 4. Universality Regained? Information & Definability in the Real Universe

What we are talking about here is a manner of material computation over higher type information: one which is framed in terms of natural language. At the classical level, natural language takes us away from computable relations very quickly. The inclusion of just one quantifier is sufficient for us to define the halting real in terms of the computable basic structure of Turing machine dynamics. This might encourage us to distrust natural language as a rather inadequate means to capture the essentially computable nature of the universe in all its variety and richness of wonders. When language does give us something incomputable in the classical sense, one might take this as evidence of the artificiality of the mathematical development, and of its enlisting of infinitary or global forms of data which are beyond our strong sense of what actually is practically computable.

This view is being weakened on a regular basis nowadays by the convergence of novel mathematical descriptions with physical and mental phenomena which the old computational framework struggles to capture.

Once again, it was Turing who was sensitive to the pulse of nature. His differential equations, based on conjectured reaction-diffusion models for the underpinnings of a range of emergent patterns in nature, were both interestingly reminiscent in syntactic form to that of equivalents of the halting real, but also surprisingly successful in releasing solutions representing the actual graphical forms observed in nature — dapplings on cows' hides, moving patterns on tropical fish, the observed variety of radiolaria formations (see Bernard Richards' article on his work with Alan Turing in the edited volume of Barry Cooper and Jan van Leeuwen on *Alan Turing — His Work and Impact*). This was just a more mathematical instance of the more general and very familiar success of natural language in capturing what we came to trust as real in the world.

Beyond the strictures of science, language was a very effective means to packaging what is real, and to reasoning with with a reassuring reliability towards a well-focused view of the world. Of course, our confidence in the reliability was often shaken by new observations and experiences.

The new paradigm being assembled tells us that we are in the presence of an increasingly pressing recognition of phenomena whose descriptions and properties deliver convergence with the familiar world of language — in a way which brings both together under an extended computational framework, involving large interactive data, and a concomitant difficulty bringing it into the sharp computational focus we expect from classical computation. Part of the uncertainty here is associated with the process of reducing or approximating the higher order information. The mathematics enveloping the higher order framework is not so easily reduced. In Newtonian times descriptions of aspects of the motions of the planets were computationally substantiated via the mathematics of the calculus. In the last century statistical analysis, game theory, differential equations and learning theory, etc., gave us descriptions of higher order relations with recognisably computational underpinnings. This expanding view of what we regard as computable — commonly taking in phenomena from economics, biology, sociology etc. which are simulatable, though not classically computable — is giving substance to what previously has seemed outside the domain of computational analysis.

Here we usefully mention the important work being carried through by the Budapest group of István Németi and Hajnal Andréka, capturing in first-order language and computational detail extreme relativistic physics. What is specially important about this work is its basis in observationally accessible physics, with theoretically mature description. We note that first-order language interpreted over large domains *can* be used to describe structures populated by what we recognise as higher type objects. What *is* limited reflects our situation as observers hosted by a particular informational context — this constituting a range of typed information constrained by our available language. Many of our examples of emergent natural relations involve a range of observational and theoretical barriers — at the quantum level, at the level of human mentality, and with the small-scale detail of emergence in biology. The work of the Budapest group on axiomatising general relativity has a distinguished anticipation in the 1920s explorations of Einstein's special theory by Hans Reichenbach (whose most famous student is no less than Hilary Putnam); while the computational aspects are thoroughly contemporary, and much dependent on our still developing understanding of the physics.

As we mentioned already, the drawing together of classical and higher order computational models was anticipated by the work of Emil Post. Post's prophetic vision of the importance of the relationship between language and computation has not just dominated the field of classical computability via 'Post's Theorem' for the arithmetical hierarchy, and via 'Post's Problem' regarding natural informational content of computational relationships at the local level: his basic sense of the key relationship carries forward productively to the newly emerging world of higher order computation as an explicitly real presence in our lives.

We still have a lot to learn about the infrastructure of higher order computation. The sense of this matches David Ruelle's comments on the mysteries of turbulence, and of the emergent relations over chaotic contexts. Other approaches to the problem have come more from those concerned to extend the classical computational models to higher levels, particularly via extensions of recursion theoretic and machine models to higher type data. Definability is just one framework within which to capture meaningful notions of higher order computability. Early work was that of Stephen Kleene, the development spread over a sequence of seminal papers.

One can find a useful review of computability-theoretic notions, related to work of such early pioneers as Kleene, Kreisel and Gerald Sacks, in John Longley's article on *Notions of Computability at Higher Types I*, from 2005. From a real-world point of view, the picture is not a pretty one. Longley develops a framework for comparing different models, and outlines the extent to which the notions have multiplied, pointing to uncertainties about what exactly are the more natural formalisations. In contrast to the classical context, the picture Longley depicts is one of a conceptual lack of robustness:

> *It is … clear that very many approaches to defining higher type computability are possible, but it is not obvious a priori whether some approaches are more sensible than others, or which approaches lead to equivalent notions of computability. In short, it is unclear in advance whether at higher types there is really just one natural notion of computability (as in ordinary recursion theory), or several, or no really natural notions at all.*

The new book by Longley and Normann on *Higher Order Computability* brings together many different strands of research in this direction. And beyond that, there are more set theoretical approaches. See Sacks' 1990 book on *Higher Recursion Theory*, or Slaman's 1981 PhD thesis on *Aspects of E-Recursion*. There is clearly much work to be done here, with a growing sense of the importance of this deep and fascinating area of investigation.

However, the natural model directly related to the science, is that derived from definability over the Turing universe. And there is clearly scope for further enrichment of this model via a drawing together of language and purely computability theoretic structure in the spirit of Emil Post.

Now is the time to return to the key concept of universality. We saw that universality at the classical level — or at least the informed view of universality at that level — is that the 'machine' has been turned into data, and we have a 'dictionary' of data-coded machines delivering in-computer universal ownership of the computing world. When we are forced by what is received as 'real' to admit higher order data, often in the form of computing 'machines' from nature — say, a rotating black hole, a body of slime mould, or a living human brain — this control through the coded list is, at the best, problematic.

There does remain a logical form to the computational environment, one which may be *described* in natural language. And looking again at the classical model, we note that it was the convergence of logical form, in the shape of the program, and the computationally significant part of the embodied machine, which gave us the universality. At higher types we have no hope of maintaining infallible control over the data, which appears in the form of found computing machinery in the real world. But we can use the logical jacket — for example, appearing in the form of the probabilities governing the incomputability of micro-analysed atomic decay — to enable us to live with nature in a productive and creative manner. We should view the moving away from the classical model, not so much as a breakdown of the 'machine as data' paradigm, but more as a loss of absolute control over our embodment of higher order data. Just as one achieves reasonably predictable outcomes through the training and riding of a horse, so computation can become a structured activity where the embodiment of the data is provided via a process which has been out-sourced to elements of the natural environment.

One might relate this to the pre-Turing experience of the word 'incomputable'. The largeness of a set can give one a very real obstacle to its being 'computed' in a natural environment, though it may be approximated. With better understanding of the activity of computation, the 'incomputability' moved up a type, and obstructed approximation. With an even higher type universe, say type-2 or more, the generic data itself becomes at best approximated via an inexact sampling.

Notice what has happened to the machine as data paradigm. It is not data which is safe, under our control. It is the logical structure we can rely on. The data has no existence disembodied. Its representational functionality inverted is its embodiment in this context. The data *is* the machine in the broadest sense. Our computer is physically changed by our import of new data, the machine protean. Things clarify when the machine is, say, the brain. For all its complexity, one can extract logic and mathematics. What is inescapable is a dimension to the data which we may only enter into via our visceral involvement with its embodiment. We retain machine as data, but as a truism without useful content — for our purposes, there is no type-defying logic in the packing and unpacking of the machine. While it is the instantiation of the logic immune embodiment that is found within the Turing universe.

The development of routes to the improvement of our grip on the slippery monster of higher order data cannot be an entirely algorithmic one. The identification of the paths through the jungle — or up the mountain — may be associated with familiar unsolvable problems from the classical level. The useful structures we use to clothe the uncertainties may need to be empirically arrived at, as are Hofstadter and Sander's analogies. Because they are embodied in physical reality, some may *have* to be empirically arrived at, if physical reality is under-determined ( such as when physical constants can have 'arbitrary' values).

In terms of the Turing universe model, we are reaching, in an inevitably ad hoc and empirical way, for definitions — ones which may be usefully retained from past — in the context of a familiar understanding of basic computable causal structure. The universality is contributed via definable relations over the Turing universe. The data — which leads to embodied local substructures of the Turing universe — is of a type beyond our computational control. It is not totally out of reach. But is forcefully calling upon sampling, approximation, statistical methods, and human cognitive resources.

There is no crisply utilisable Church-Turing thesis at the higher levels we encounter in the physical world. The brain is awesomely effective, applying its corpus callosal coordination between classical and real-world-engaged higher order computation. We now achieve a mature context in which to re-focus on Turing's recognition that (quoting again from his Lecture to the *London Mathematical Society* on February 20, 1947):

> *… if a machine is expected to be infallible, it cannot also be intelligent. There are several theorems which say almost exactly that.*

Turing's acceptance of the way in which humans can enrich algorithmic 'thinking', on which his vision of cooperation between computer and human intelligence is based, is related to the mathematics in his popular piece on *Solvable and Unsolvable Problems* (*Penguin Science News* 31, 1954, final paragraph of p. 23):

> *The results which have been described in this article are mainly of a negative character, setting certain bounds to what we can hope to achieve purely by reasoning. These, and some other results of mathematical logic may be regarded as going some way towards a demonstration, within mathematics itself, of the inadequacy of 'reason' unsupported by common sense.*

This aspect of Turing's thinking is not so much remembered nowadays, but takes on a special significance in a contemporary setting. This is the point at which, giving due respect to the mathematics of computation, higher order data and definability, we ask what sort of specific benefits we might expect from our computational model of causality and its inductive type structure.

## 5. Some Rewards

What we have attempted to describe is a world in which information obeys basic rules of typing, with the scaling of different information levels being characterised by mathematically simple amalgamations. We have pointed to connections between the mathematical structures and the more descriptive analysis of information of Floridi and others. We have discussed the rooting of computability in the relativity of the observer to contextual information, its relationship to what we view as 'causality', and the epistemological and causal constraints modelled via the Turing universe. We have persistently stressed the universality of the extended Turing model, with its close relationship to the most fundamental interactivity we encounter within our observational ambit. This elevates the analysis to a modelling exercise rather than one driven by a mining of metaphors or analogies. The mathematical transitions between types of information as embodied in the physical and mental universe - aspects of the same world - have been packaged within various more or less familiar frameworks. Most evocative in relation to the philosophy, arts and humanities is the neglected notion of 'definability', and it is this we have exercised productively over the basic Turing model of computable causality. The motivation has been to give mathematical substance to the empirically arising natural laws (in particular), and give an overarching context for the emergence of real relations from their more basic physical contexts. What is a particularly striking aspect of this developing picture is the drawing together of ostensibly conflicting activities, as pursued by communities representing C.P. Snow's 'two cultures'. The relationships between language, definability, computation, biology, physics, mentality and machine begin to emerge as aspects of a general framework within which one can locate the particularities of diverse causal contexts — while at the same time, arriving at a less linear and mechanical concept of causality.

We finish with an indicative check-list of questions and issues which may be clarified via the accelerated computational understanding that we have described above.

● **The fragmentation of science** we discussed in some detail above. The key ingredients we identified were global outcomes and corresponding descriptions, involving the discerning of phase transitions and emergent entities and relations; and the intervention of mechanisms for type reduction restoring an apparent informational 'level playing field', accessible at our own level of existence, mediated by a residue of computability theoretic impediments to reduction between LoAs. Returning to the relationship between physics and biology, we are happy to accept that despite the causal connections between particle physics and the study of living organisms, the corresponding disciplines are based on quite different basic individuals and natural laws, and there is no feasible and informative reduction of the higher level to the more basic one.

Such considerations feed into issues related to some of Lee Smolin's 'Great Problems', from his 2006 book on *The Trouble with Physics*:

● **Combining of relativity and quantum theory into a single theory** that can claim to be the complete theory of nature. Smolin characterises this (*The Trouble with Physics*, p.5) as being 'the problem of quantum gravity'. The Turing model currently supports the dichotomy between a 'low level' or 'local' structure with a much sparser level of uniquely definable relations than one encounters at higher levels of the structure, where the richer context translates into more coding resources. Simply put, the reason for this situation is that the higher one ascends the structure from computationally simpler to more informative information, the more data one possesses with which to describe structure.

One should remind the reader that the focus is on *models* of 'reality'. The argument of this article is that different empirically consistent models may point to incompletely defined underlying physical characteristics, using 'defined' in a sense appropriate to the current level of knowledge of the relationship between relevant relational structure and language. Of course, the validation of modelling itself cannot be proved, adding a level of uncertainty to the history within which the currently appropriate model is permitted to function — one could mention the Duhem-Quine thesis. One should be careful not to allow type-theoretic notions of 'higher' or 'lower' to entail some expectation or implication of causal 'dominance'. This may fit with a desire to linearise conceptual causality, but is not appropriate to a context in which one regards the structure as expressing an organic interdependence which is far from linear, except in the sense that the observer may get impressions of dominance according to the pertaining location.

The conclusion which beckons is that there is no *useful* theory unifying relativity and quantum theory — the fragmentation of nature is built in, via its mathematical underpinnings. Useful is the key word here. One might extend our description of the relationship between quantum and classical levels, via descriptions which do not improve practical transfer of predictive information across boundaries. But this does not preclude the discovery of an enhanced theoretical framework bringing great insights to the physics. We might add as a footnote that even the question of *consistency* of quantum theory with special relativity and interaction presents serious problems. The mathematician Arthur Jaffe, with J. Glimm and other collaborators, succeeded in solving this problem in space-time of less than four dimensions in a series of papers (see Jaffe's *Quantum Theory and Relativity* from 2008), but it is still open for higher dimensions.

● **Finding a realistic foundations of quantum mechanics** is essentially Smolin's second question. Under this heading we can group a number of individual though related issues. The mysteries which haunt the realist include: non-locality; ontological ambiguity and entanglement; and decoherence arising from collapse of the wave function. But — relative to an acceptance of what we can observe, and of the possibility that there are no foundations — all can be made sense of under the heading of *immanence*. We start from an assumption that the higher order relations on a structure are those arising from its more basic properties. That is, we assume that the material universe has to do everything for itself. It has nothing it does not define.

What is required is the taking on-board of definability, and invariance under automorphisms, as the work-horses of immanent development of our universe. And key to these is the nature of the automorphism group of the Turing universe.

The reader might ask: Would not *any* mathematical structure with non-trivial automorphisms furnish an equally good analogy for the physical situation in question? Of course, by isolating particular contexts one may indeed pin down mathematical structure with a good match between mathematical properties and empirical evidence. What is more fundamentally satisfying about the Turing model is its independence of arbitrary assumptions, the context-free nature of its functionality, in a way that even string theory has not managed to deliver on, in particular in its failure to pin down the geometry of space-time. It is the mathematically solid basis of the Turing universe in the computational content we observe, and nothing else, which promises to deliver structure across a swathe of scientific fields. In doing this, we mathematically access the transitional passage between scientific disciplines via the computational/definitional infrastructure of the typing of information, and observe the restoration to scientific feasibility at different levels via those same emergent/ definitional relations that create the ontological elements, and precisely capturing the nature of the actual mathematics used to predict from tried and tested natural laws, crucial to our context-free approach from below. The down-side is that this solidity of foundational content is received at high cost to the specificity deliverable at this point in time. Even basic questions concerning the automorphism structure of the Turing universe remain unresolved. In a sense, this is promising for those who have respect for the complexity of reality as we experience it.

If the whole structure can be mapped onto itself, associating a particle with two different locations in space-time , then the particle can but exist in two different locations. The real universe has no means to prevent this. If entanglement changes the role of an aspect of the universe in relation to the automorphisms permitted, so does the ontology and the associated epistemological status. Observation can do this to one of our particles. Of course, the lack of a comprehensive description of the particle will have to be shared by whatever form the particle takes. There is nothing odd about the two slit experiment. Or about the dual existence of a photon as a wave or as a particle. And the interference pattern is just evidence that a particle as two entities can define more than it can solo. And as for decoherence, many-worlds and the multiverse? It is not just that the whole scenario becomes redundant, it is that the acceptance of such easily achieved permutations of reality is mathematically naive. A small change in a complex structure with a high degree of interactivity can have a massive global effect. And the mathematics does have to encompass such modifications in reality.

Realistically, one can only discuss the case of the 'branching multiverse' arising from its many-worlds origins. The decoherence is triggered by the intervention of a measurement, observation, or similar. This intervention has to be sufficient to induce an inconsistent reality entailing non-trivially interesting decoherent branches. The mathematical experience with definability and invariance over interesting structures leads us to envisage an extended relational structure encompassing such far-reaching mathematical consequences that *either* the branching becomes unlikely to be comprised of sufficiently interesting structures as to give the branching any great functional role in its contribution to an interpretation with a claim to realism; *or*, more persuasively, the alternative branches merely disappear in their entirety, leaving behind a common-or-garden imposition of extended Turing definability.

The mathematics does have to encompass such modifications in reality? More to the point, its character is established by the nature of its automorphisms, So — there is a qualitatively different apparent breakdown in computability of natural laws at the quantum level — the measurement problem challenges us to explain how certain quantum mechanical probabilities are converted into a well-defined outcome following a measurement. In the absence of a plausible explanation, one is denied a computable prediction. The physical significance of the Turing model depends upon its capacity for explaining what is happening here. If the phenomenon is not composite, it does need to be related in a clear way to a Turing universe designed to model computable causal structure. We look more closely at definability and invariance.

● **Do we need the multiverse?** Let us first look at the relationship between automorphisms and many-worlds. When one says "I tossed a coin and it came down heads, maybe that means there is a parallel universe where I tossed the coin and it came down tails", one is actually predicating a large degree of correspondence between the two parallel universes. The assumption that you exist in the two universes puts a huge degree of constraint on the possible differences — but nevertheless, some relatively minor aspect of our universe has been rearranged in the parallel one. There are then different ways of relating this to the mathematical concept of an automorphism.

One could say that the two parallel worlds are actually isomorphic, but that the structure was not able to define the outcome of the coin toss. So it and its consequences appear differently in the two worlds. Or one could say that what has happened is that the worlds are not isomorphic, that actually we were able to change quite a lot, without the parallel universe looking very different, and that it was these fundamental but hidden differences which forces the worlds to be separate and not superimposed, quantum fashion. The second view is more consistent with the view of quantum ambiguity displaying a failure of definability. The suggestion here being that the observed existence of a particle (or cat!) in two different states at the same time merely exhibits an automorphism of our universe under which the classical level is rigid (just as the Turing universe displays rigidity above $0''$, the Turing degree of the Halting Problem relativised to itself, see Slaman [55]) but under which the sparseness of defining structure at the more basic quantum level enables the automorphism to re-represent our universe, with everything at our level intact, but with the particle in simultaneously different states down at the quantum level. And since our classical world has no need to decohere these different possibilities into parallel universes, we live in a world with the automorphic versions superimposed. But when we make an observation, we establish a link between the undefined state of the particle and the classical level of reality, which destroys the relevance of the automorphism.

To believe that we now get parallel universes in which the alternative states are preserved, one now needs to decide how much else one is going to change about our universe to enable the state of the particle destroyed as a possibility to survive in the parallel universe — and what weird and wonderful things one must accommodate in order to make that feasible. It is hard at this point to discard the benefits brought by a little mathematical sophistication, with Occam's razor on our side. As George Ellis puts it in *The Unique Nature of Cosmology* in 2003:

> *The issue of what is to be regarded as an ensemble of 'all possible' universes is unclear, it can be manipulated to produce any result you want … The argument that this infinite ensemble actually exists can be claimed to have a certain explanatory economy (Tegmark 1993), although others would claim that Occam's razor has been completely abandoned in favour of a profligate excess of existential multiplicity, extravagantly hypothesized in order to explain the one universe that we do know exists.*

Quantum ambiguity as a failure of definability is a far more palatable alternative than the invention of new worlds of which we have no evidence or scientific understanding.

● **Explain how the values of the free constants in the standard model of particle physics are chosen in nature.** This is Question 4 from Smolin, echoing Einstein in his *Autobiographical Notes*, to be found in *Albert Einstein: Philosopher-Scientist*, 1969, p.63):

> *… I would like to state a theorem which at present can not be based upon anything more than upon a faith in the simplicity, i.e. intelligibility, of nature … nature is so constituted that it is possible logically to lay down such strongly determined laws that within these laws only rationally completely determined constants occur (not constants, therefore, whose numerical value could be changed without destroying the theory) …*

Again, the drawing together of a global picture of our universe with the basic mathematical model of the Turing universe flows from the correspondence between emergent phenomena and Turing definable relations. This gives us a framework within which to explain the particular forms of the physical constants and natural laws familiar to us from the standard model science currently provides. It goes some way towards substantiating Penrose's 'strong determinism' (see his *Quantum physics and conscious thought,* from *Quantum Implications: Essays in honour of David Bohm*, 1987, pp. 106-107*)*, according to which:

> *... all the complication, variety and apparent randomness that we see all about us, as well as the precise physical laws, are all exact and unambiguous consequences of one single coherent mathematical structure.*

Included in this picture must be the problem of the geometry of space-time, which Smolin also spotlights in his dismissal of string theory and his advocacy of other approaches, such as that of quantum gravity.

Smolin himself has travelled some distance since his *The Life of the Cosmos* — his new viewpoint appears fully developed in his (joint with Roberto Mangabeira Unger), *The Singular Universe and the Reality of Time: A Proposal in Natural Philosophy.* Here there is a continuity of vision of globally determined natural constants, geometry and laws, but within a uniquely determined universe — the 'sufficient reason' is self-applied by the universe, in keeping with the view adopted here.

We leave the physics, and Smolin's other questions, and look very briefly at some further real world puzzles. These share with those from physics the characteristic of being of great interest to large communities of researchers, but largely approached in very specific contexts, but with little clear sense of the basics of what they are dealing with. The knowledge of the general concepts and technical frameworks flowing from logic is of a generally low level. We are in the age of hugely energetic and scientifically prolific academic ants — ants who build but have not been trained to appreciate the informational content of the intricate and formal patterns enlisted by their own Stakhanovite achievements. A particular example is that of network theory, the specificity of the research both unavoidable for theoretical reasons, while providing the key to our interest of the area. The modelling and frequent surprises are their own reward. But for those with a need for activity augmented by consciousness, there remains the question of:

● **What exactly is emergence, and in what sense does it exist?** Ronald Arkin's observations on emergence in *Behaviour-Based Robotics* are not easily dismissed:

> *Emergence is often invoked in an almost mystical sense regarding the capabilities of behavior-based systems. Emergent behavior implies a holistic capability where the sum is considerably greater than its parts. It is true that what occurs in a behavior-based system is often a surprise to the system's designer, but does the surprise come because of a shortcoming of the analysis of the constituent behavioral building blocks and their coordination, or because of something else?*

In a 1999 *Artificial Life* article (*Design, observation, surprise! A test of emergence*, pp. 225–239) Ronald, Sipper and Capcarrère clarified the common view of emergence, formulating a 'test' for emergence modelled on the Turing test for intelligence — but without really answering the question. They nicely picked out the need for different languages to describe the basic causality and the emergent phenomena, respectively. But they fell back on the old criterion of 'surprise', questioned by Arkin, for an outcome to be emergent. Our characterisation in terms of definability, and the nature of the language needed for the definition, does the job without such appeals to subjective judgements. Of course, there are situations where we might struggle for a definition to analyse, where we might have doubts about the level of definition actually required, or the logical complexity identified may not deliver much surprise to the typical observer. However, one may expect a good correlation between the observer-led expectation of an outcome arising from an essentially global collation of the basic causal structure, and the logical character of the description of it.

A particularly well-trodden and challenging area of application of the emergence-derived conceptual framework relates to the brain and human mentality. We look briefly at some commonly considered aspects related to intelligence, consciousness and free-will.

● **How exactly does mentality supervene on brain activity? And is there more than subjective content to our recognition of emergence at work?** This is clearly a large topic, and this is not the place to focus in depth on the issues. We just comment briefly about the appropriateness of mathematics in such an informationally messy and mysterious context, and about the clarifying benefits of having a basic model to host one's thinking. There can be no more persuasive example of the machine as data, and of the regard due to informational structure. There are kinds of thinking that are popularly held to be directly ruled by rationality, and which might more easily be carried out by a machine than others which we might ascribe to 'common sense' or 'invention'. In the popular mind, mathematical deduction probably comes top amongst such candidates for thinking reducible to pure reason.

This, of course, was the basis for David Hilbert's programme aimed at replacing mathematicians — or at least mathematics — with grand formal theories following set rules of deduction. Such thinking may well have been what drove Alan Turing to try to transcend Gödel's incompleteness theorem via a hierarchy of constructively defined extensions of incomplete theories, indexed into the transfinite by Kleene's recursive ordinals. So it is to mathematicians we particularly turn — to Henri Poincaré in his 1908 lecture at the *Société de Psychologie* in Paris, and to Jacques Hadamard, carrying the baton into the post-war years with his 1945 "*The Psychology of Invention in the Mathematical Field*". Hadamard's account of Poincaré's celebrated moment of mathematical insight on stepping onto a bus is worth repeating:

> *At first Poincaré attacked [a problem] vainly for a fortnight, attempting to prove there could not be any such function … [quoting Poincaré]:*
>
> *"Having reached Coutances, we entered an omnibus to go some place or other. At the moment when I put my foot on the step, the idea came to me, without anything in my former thoughts seeming to have paved the way for it … I did not verify the idea … I went on with a conversation already commenced, but I felt a perfect certainty.*
>
> *On my return to Caen, for conscience sake, I verified the result at my leisure."*

A striking part of the story is the "certainty" — the sense that there was a memetic quality to the apprehension, consistent with the existence of a representational handle on the form of the proof … something close to an algorithm for generating it from the key ingredient. Of course, 'certainty' may have its origins in a sampled factual context, where the selection of information turns out to be misleading and unrepresentative in some key way, leading to what is later realised to be a wrong conclusion. The importance of the "certainty" here is not in relation to the identification of 'truth', but as a recognisable symptom of a certain logical character to the mental process, a memetic definability, a complete conceptual structure, where the association with something valid in the longer term is characteristic of a fully realised and characterised logic, something familiar and unmistakable to the working mathematician.

The Turing model is there in the Poincaré story. The non-algorithmic search for the mental formation fitting the hard-won accumulation of mathematical detail, some relevant, some not. We are struck by the fact that the formation was a *certain* fit, even though not delivered in precise detail, the detail to be put in place at a later date. And — we see — the diversion presented by the practicalities of the journey just enough to free up the mind to carry out the sort of higher order processing of information needed. But:

● **What evidence is there of different kinds of thinking, both higher order and algorithmic? What can we say about higher order mentality, its representation, and downward causation?** Back with the realities of the physical world, we do see evidence of the partnership between different kinds of thinking being mirrored within the structure of the human brain. In the Turing model — not just the pure 1936 model underpinning the informational flat-earth, but the interactive oracle machine based model — there occur the higher order definitional accoutrements which both correspond to a suitably embodied Poincare phenomenon. And it is this separation of computational content which we see vividly reflected via the neuroscience. In the model we harvest the higher order outcomes when not busy with the underlying machinery, and vice versa. Notice that the mathematics points to the possibility of these two worlds converging on overlapping informational domains. Consequent benefit and coherence requires suitably embodied coordinating mechanisms, and this too we see.

The brain is neatly structured into physically similar right and left hemispheres, but with the functional characters which feature so prominently in popular takes on human thinking. The popular thinking is reflected in more considered takes on the subject, such as this taken from the 2009 book *The Master and his Emissary: The Divided Brain and the Making of the Western World* by Iain McGilchrist:

> *The world of the left hemisphere, dependent on denotative language and abstraction, yields clarity and power to manipulate things that are known, fixed, static, isolated, decontextualised, explicit, disembodied, general in nature, but ultimately lifeless. The right hemisphere by contrast, yields a world of individual, changing, evolving, interconnected, implicit, incarnate, living beings within the context of the lived world, but in the nature of things never fully graspable, always imperfectly known — and to this world it exists in a relationship of care. The knowledge that is mediated by the left hemisphere is knowledge within a closed system. It has the advantage of perfection, but such perfection is bought ultimately at the price of emptiness, of self-reference. It can mediate knowledge only in terms of a mechanical rearrangement of other things already known. It can never really 'break out' to know anything new, because its knowledge is of its own representations only. Where the thing itself is present to the right hemisphere, it is only 're-presented' by the left hemisphere, now become an idea of a thing. Where the right hemisphere is conscious of the Other, whatever it may be, the left hemisphere's consciousness is of itself.*

It is the corpus callosum, a feature of the brain architecture of placental mammals in general, that connects and mediates the functionality of the separate hemispheres. McGilchrist describes (pp. 18-19) exactly what one might expect from our computational model:

> *... the evidence is that the primary effect of callosal transmission is to produce functional inhibition.*
>
> *... it turns out that the evolution both of brain size and of hemisphere asymmetry went hand in hand with a reduction in interhemispheric connectivity. And, in the ultimate case of the modern human brain, its twin hemispheres have been characterised as two autonomous systems.*
>
> *So is there actually some purpose in the division of neuronal, and therefore, mental processes? If so, what could that be?*

The corpus callosum, one assumes, plays a key role in the interchange of the fruits of different levels of thinking, incorporating representational devices for traversing the frontiers between different typed classes of information. These two quotations of McGilchrist provide us with a conceptual bridge between two quite different areas of experience and analysis, with dramatically convergent logical characters. From McGilchrist we get a physical counterpart to our cultural experience of different types of thinking (with the dual use of 'type' being intentional). From the mathematics we get a still inadequately understood and explored convergence of different levels — or more technically speaking now — *types* of computation with both complementary and disjunctive relationships.

The mathematics confirms intuitions about the relationship in the brain, with the crucial role identified for the corpus callosum. And the confusions and acres of descriptive speculation surrounding the physical and cultural engagement seems to cry out for the clarifying vision delivered by the mathematical analysis of computational and informational structure. What Newton and Einstein brought to time and space is overdue now in the context of information and its causal/computational structure.

We can see the embodied mathematics of representation emerging from Antonio Damasio's description (p.170 of *The Feeling Of What Happens*) as he describes how:

> *… both organism and object are mapped as neural patterns, in first-order maps; all of these neural patterns can become images …The sensorimotor maps pertaining to the object cause changes in the maps pertaining to the organism … [These] changes … can be re-represented in yet other maps (second-order maps) which thus represent the relationship of object and organism …The neural patterns transiently formed in second-order maps can become mental images, no less so than the neural patterns in first-order maps.*

The picture is one of re-representation of neural patterns formed across some region of the brain, in such a way that they can have a computational relevance in forming new patterns. So far missing are the specifics of the computational host. These do exist, consolidating the relevance of the Turing universe as their host. The early proposals for connectionist computational models of the synaptic brain structure go back of course to the McCulloch and Pitts 1943 article on *A logical calculus of the ideas immanent in nervous activity* — though Turing had an independently devised version (see Christoff Teuscher's 2002 book on *Turing's Connectionism*).

● **What does the Turing model add to our understanding of both neural nets and their appropriateness as a model? How does the unconscious fit into this picture?** Neural nets fit well within the Turing model, though the higher order outcomes and the associated analysis of definability and representations of higher order outcomes for re-integration into the basic computable causality is not explicitly part of the model. Without this, we see something common with new computational paradigms, in that they connect up with higher order computational outcomes without being specifically adapted to make use of them. There are, as we have already seen, fundamental differences between computation at the classical level, and what happens practically in relation to computation over higher type information. And that even the brain has evolved to separate out these different computational modes, accompanied by sophisticated means to the coordination and separation of these modes.

One can generally group, intuitively if not formally, the associated candidates for embodying classically incomputable outcomes along with the halting real. Given that many of these 'new paradigms' are based on familiar embodied computation in nature, it is no surprise that the non-classical outcomes one observes *are* observed at all, with far more embodiment than the purely abstract haling real. Here is one credible example of speculation in this direction, from Paul Smolensky in 1988, before he moved to a more exclusive focus on linguistics (from *On the proper treatment of connectionism*):

> *There is a reasonable chance that connectionist models will lead to the development of new somewhat-general-purpose self-programming, massively parallel analog computers, and a new theory of analog parallel computation: they may possibly even challenge the strong construal of Church's Thesis as the claim that the class of well-defined computations is exhausted by those of Turing machines.*

What seems clear, both from our intimate experience of brains and their associated mentality, and from what we have in the way of informative computational models, is that the complexity of interaction between different kinds of thinking, modeled by computation over differently typed information, is key to what is special about human mentality. The nonlocal computation, equipped with representational infrastructure and looping between types is what is hard to reproduce structurally. As it is, we have Rodney Brooks proclaiming in Nature in 2001 that:

> *… neither AI nor Alife has produced artifacts that could be confused with a living organism for more than an instant.*

It would be hard to disagree.

We can argue about the reality or otherwise of consciousness, and the subjective feeling that it is far from illusory. What is the case is that the mathematics, complete with typed structure of information and computational richness of relationships, is capable of providing a fitting host for an emergent level of information which may be supervenient on the physical brain, but is not subsumed by the basic properties of it. The evidence of the representational activity embodied in a manner identified by Damasio and modelled via the Turing universe and its hosting of higher type outcomes captured in definitional packages is very persuasive.

Not all this embodied computational activity is reported back. The structural connectivity has many nuances, and as one practicing neuroscientist commented after a recent talk "the more one learns about the brain, the more one understands how much more there is to learn". Certain thoughts may need locating or hiding, like the hidden chaos of domestic cupboards and drawers. Here we have Daniel Dennett, known reductively in the popular mind for a reductionist attitude to the conscious and unconscious, referring in his 1991 *Consciousness Explained* (p. 308) to unconscious 'thoughts' as a higher order phenomena:

> *Unconscious thoughts are, for instance, unconscious perceptual events, or episodic activations of beliefs, that occur naturally — that must occur — in the course of normal behavior control. Suppose you tip over your coffee cup on your desk. In a flash, you jump up from the chair, narrowly avoiding the coffee that drips over the edge. You were not conscious of thinking that the desk top would not absorb the coffee, or that coffee, a liquid obeying the law of gravity, would spill over the edge, but such unconscious thoughts must have occurred — for had the cup contained table salt, or the desk being covered with a towel, you would not have leaped up. Of all your beliefs — about coffee, about democracy, about baseball, about the price of tea in China — these and a few others were immediately relevant to your circumstances. If we were to cite them in an explanation of why you leaped up, they must have been momentarily accessed or activated or in some way tapped for a contribution to your behavior, but of course this happened unconsciously.*

All this makes daunting reading for those hoping to construct intelligent machines. So:

● **What are the prospects for artificial intelligence? And what are the pointers to an appropriate relationship between computer and human intelligence?** Those working in AI might be encouraged by the existence of a model which roots intelligent thought in some sort of embodied context, however resistant it may be to simulation. It may be some comfort to know: Even if the brain does represent and utilise more of its higher order outcomes than the termite in her cathedral does, the human brain taking possession not just of the 'syntax' but of some accessible level of its accompanying 'semantics' in the form of higher type information, all this *is* embodied. The problem for AI is not so much the theoretical barriers to embodying intelligence — the problems are practical rather than absolute. We are back where the 1940s builders of the universal computer were, seeking to build an artifact, predicted by theory, but whose embodiment was a whole new adventure.

The brain is itself emergent, computed by a world deploying huge informational resources, in both space and time. The successes of today's artificial intelligence are real, if far short of Turing's Hanslope Park ambition of "building a brain". When Marvin Minsky at Boston University in May 2003 declared that:

*AI has been brain-dead since the 1970s*

the hyperbole drew attention to an important observation. AI was being pursued without any real understanding of the fundamentals of the activity. Did this matter? Maybe we *would* have got the stored program computer without the help of Turing. It is true that even today many very able people think they understand the universal Turing machine and its theoretical underpinning of the modern computer, but really don't. Universality is indeed a short step from programmability — for who knows how to *program in* the universality. The universality facility is basic — as potentially basic and memetic and ingenious in conception as is the wheel. It does not happen by accident or incrementally, even though an ad hoc process may host hidden logical complexity. The transition from rolling tree trunk to wheel is a major intellectual jump. As is Turing universality, and its dependence on a representational sophistication capable of hosting the 'program as data' paradigm.

What has always been true is that enhanced *consciousness* of achievements and limitations is the route beyond the termite cathedral. Or, in this case, the *Turing's Cathedral* of George Dyson's popular book. One corollary of all this is a validation of Turing's view of the necessity of cooperation between machine and humans, of course. Before a nod in the direction of the free will debate, we allow some *very* brief comments on the causal relationship between brain and mentality: the old problem passed down to us from the time of Rene Descartes.

● **How can the computational model host diverse levels of causality coherently without over-causation? What does the solution to this problem tell us about the non-linearity of causality?** Susan Blackmore is an articulate representative of those who would avoid difficulties via a retreat to informational flatness (page 220 of *The Meme Machine*, Oxford University Press, 1999):

> *Dualism is tempting but false. For a start no such separate [thinking] stuff can be found. If it could be found it would become part of the physical world and so*

*not be a separate stuff at all. On the other hand if it cannot, in principle, be found by any physical measures then it is impossible to see how it could do its job of controlling the brain. How would immaterial mind and material body interact? Like Descartes' 'thinking stuff', souls, spirits and other self-like entities seem powerless to do what is demanded of them.*

Over many years and many books, the philosopher Jaegwon Kim has frequently revisited such questions as: How can mentality have a computational role in a world that is fundamentally physical? And what about 'overdetermination', the problem of phenomena having both mental and physical causes? In *Physicalism, or Something Near Enough* (Princeton, 2005) Kim phrases the dilemma as:

> *… the problem of mental causation is solvable only if mentality is physically reducible; however, phenomenal consciousness resists physical reduction, putting its causal efficacy in peril.*

We have already entered a different world, the mathematics helping us leave behind the linear causality of everyday practicality. What the extended Turing model gives us is an organic whole in which 'causality' is a consequence of informational dynamics beyond the linearity of popular thinking. No wonder our schooldays, where we tried to mould the fourteen reasons for the start of the Great War of 1914-18 into some coherent understanding of whose fault it was, were so frustrating. No wonder it is so hard to understand why great civilisations have given rise to a history littered with such monstrous crimes against humanity. There is 'causality' which can be decoded at different levels, appropriate to the observer's vantage point. The causality is not illusory. It derives its power from mathematical relationships yielding up a degree of predictive power. The organic whole grants us this, and may take it away. Or it may just hide the mathematical information behind a turbulent confusion of computational interactivity. We ourselves may be part of this turbulence.

● **What part does free will play?** As for free will, the attendant problems of definition are similar to those we face in defining randomness. Most discussions bypass the difficulties of defining what it means for an embodied dynamic to be somehow selectively embedded in its causal context. Just as there is no meaningful definition of absolute randomness, we expect the notion of 'free' will to be equally unformulatable. What we can describe is a place for computational activity associated with the reality of a 'machine' such as the brain, where the computational structure hosts a level of outcome via a process only computable at a level not accessible to the computational resources of the decision framework. Typically, one might conceive of an emergent decision only achieved via considerable conscious (that is, reported) thought processes, but arising out of a wide reaching ('global') context in a — 'surprising', we allow ourselves — way. It 'feels', by what is reportable of this higher order computation, to be surprisingly related to our (reported) sense of identity. It does feel like free will. And this is an integral part of our engagement in the process.

There are many useful writings on the topic of free will, of course. We just mention one which came our way very recently, by way of a Turing Centenary themed issue of the Philosophical Transactions of the Royal Society, series A, where Seth Lloyd writes engagingly on *A Turing test for free will*.

What we are left with is a sense of a more coherent and exciting universe, of which the computer does little more than scratch the surface. It is a world in which information is a discipline problem, but an equal partner with logical structure. It is one in which there is room for the categorical quantum theorist to extract beautiful analogies between different contexts; or for the application of large computational resources in service of the taxonomy of myriad emergent wonders of informational specificity. The taming of data, versus visceral engagement? — merely two sides of the same coin.